\documentclass{amsart}
\usepackage{amsmath,amssymb,amsfonts,graphicx}
\usepackage[all]{xy}

\newcommand{\zln}{\mu_{l^N}}
\newcommand{\p}{\mathfrak{p}}\newcommand{\q}{\mathfrak{q}} 
\newcommand{\F}{\mathbb{F}}
\newcommand{\Km}{{K^\times}}
\newcommand{\Kvs}{K_v^\times}
\newcommand{\rKS}{\sqrt[l^N]{K_S}}
\newcommand{\rmKS}{\sqrt[m]{K_S}}
\newcommand{\rKSl}{\sqrt[l^{N}]{K_S}}

\newcommand{\N}{R}
\newcommand{\M}{\mathcal M}

\newcommand{\Omegae}{\Omega_{l}}

\newcommand{\Q}{\mathbb{Q}}

\newcommand{\Z}{\mathbb{Z}}

\newcommand{\es}[3]{1 \rightarrow #1  \rightarrow #2  \rightarrow #3  \rightarrow 1}
\newcommand{\les}[3]{1 \rightarrow #1  \rightarrow #2  \rightarrow #3}
\newcommand{\Frob}[2]{Frob({#2},#1)}

\newcommand{\cep}[6]{ \begin{xymatrix}{&&&G_K \ar@{-->}[dl]_{#4} \ar [d]_{#5}\\
1 \ar[r] & #1 \ar[r] & #2 \ar[r]_{#6} & #3 \ar[r]
& 1.}\end{xymatrix}}

\newcommand{\pcep}[6]{ \begin{xymatrix}{&&&D_\p \ar@{-->}[dl]_{#4} \ar [d]_{#5}\\
1 \ar[r] & #1 \ar[r] & #2 \ar[r]_{#6} & #3 \ar[r]
& 1.}\end{xymatrix}}

\renewcommand{\proof}{{\bf{\emph{Proof:\\}}}}
\newcommand{\done}{\hfill $\blacksquare$ \newline}

\newcommand{\hs}{\hspace{.2 cm}}

\numberwithin{equation}{subsection}

\newtheorem{thm}{Theorem}[section]
\newtheorem{lem}[thm]{Lemma}
\newtheorem{prop}[thm]{Proposition}
\newtheorem{cor}[thm]{Corollary}

\newtheorem{definition-lemma}[thm]{Definition-Lemma}
\theoremstyle{definition}
\newtheorem{defn}[thm]{Definition}

\theoremstyle{remark}
\newtheorem{rem}[thm]{Remark}

\newtheorem{claim}[thm]{Claim}

\title{Minimal Ramification in Nilpotent Extensions}
\author{Nadya Markin}\thanks{Research supported by the Claude
Shannon Institute, Science Foundation Ireland Grant 06/MI/006}
\email{nadyaomarkin@gmail.com}
\author{Stephen V. Ullom}
\email{ullom@math.uiuc.edu} 
\address{Nadya Markin \\CASL UCD, Belfield, Dublin, Ireland}
\address{Stephen Ullom \\1409 West Green St
University of Illinois
Urbana, Il 61801, USA}
\begin{document}

\begin{abstract}

Let $G$ be a finite nilpotent group and $K$ a number field with torsion
relatively prime to the order of $G$. By a sequence of central group
extensions with cyclic kernel we obtain an upper bound for the minimum
number of prime ideals of $K$ ramified in a Galois extension of $K$ with Galois group isomorphic to 
$G$. This sharpens and extends results of Geyer and Jarden and of Plans. Also
we confirm Boston's conjecture on the minimum number of ramified
primes for a family of central extensions by the Schur multiplicator.

\end{abstract}

\maketitle

\section{Introduction}
Given a number field $K$ and a finite group  $G$ an important problem is to find a Galois extension $L$ of
 $K$ such that its Galois group $Gal(L/K)$ is isomorphic to $G$. Scholz and Reichardt (see  Serre \cite{SerreTGT} for a modern account) proved independently that any $l$-group $G$, $l$ an odd prime, occurs as the Galois group of an extension of the rationals. Shafarevic \cite{Sha} has shown  for any solvable group $G$ and  number field $K$ that there exists a Galois extension $L/K$ with $G \cong Gal(L/K)$. 
In this paper we ask for given $K$ and nilpotent $G$, what is the minimum number $minram_K(G)$ of prime ideals of $K$  ramified in $L$ as
 $L$ runs over extensions of $K$ with $Gal(L/K) \cong G$? We rephrase the question for $l$-groups $G$: 
  For a given finite set $S$ of prime ideals of $K$, $K(l,S)$ denotes the maximal $l$-extension    of $K$ which is unramified outside $S$. How large must $S$ be so that $G$ is isomorphic to a quotient group of $Gal(K(l,S)/K)$ for some $S$? 
   
   One knows $minram_{\Q}(G)\leq n$ if $G$ is an $l$-group of order $l^n$, $l \not = 2$, cf. \cite{SerreTGT}. If $G$ is an abelian group, an application of class field theory (Theorem \ref{AbG}) shows $minram_K(G) \leq d(G):=$ minimum number of generators of $G$. 
   In fact for the case when $K = \Q$, Boston's conjecture \cite{BM} stated below implies that $minram_{\Q}(G) \leq d(G)$ for all finite groups $G$.

Suppose $G$ is a nilpotent group and the field $K$ satisfies for each prime $l$ dividing the order $|G|$ of $G$ conditions 

\begin{enumerate}
\item $K$ does not contain a primitive $l$-th root of unity $\zeta_l$
\item $K$ has no ideal classes of order $l^2$,
\end{enumerate}

then Theorem \ref{NilThm} states
$$ minram_K(G) \leq \sum_{i \geq 1} d(G_i/G_{i+1}) +t(K).$$

Here $\{G_i\}$ is the lower central series of $G$, $G_1 = G$, $G_{i+1}=[G_i, G]$ and $t(K)$ is a constant depending only on $K$. 
This extends Plans' \cite{Plans}  result on $minram_{\Q}(G)$ to all number fields $K$ satisfying $(1), (2)$ above. Secondly, Geyer and Jarden  \cite{GJ} obtain the bound $minram_K(G) \leq n+ t(K)$, where the $l$-group $G$ has order $l^n$ and $\zeta_l \notin K$. We obtain the improved bound by considering central embedding problems with a cyclic kernel, not just kernel of prime order as in \cite{GJ}. 
Note that without condition $(2)$, the methods of Section $8$ still generalize the results of Geyer and Jarden \cite{GJ} to nilpotent groups, giving a weaker bound for a nilpotent group $G$ of order $\prod_{l\mid |G|}l^{n_l}$, namely 
$$ minram_K(G) \leq \max_{l \mid |G|} \{n_l\} +t(K).$$ 

We generalize Geyer and Jarden's definition of an exceptional set $T$ of primes to the prime power setting in Section $4$; this provides the technical tool for constructing idele class characters with strictly controlled ramification.

The realization of $l$-groups is carried out in three steps similarly to \cite{GJ}, \cite{SerreTGT}, \cite{Plans}: the first step involves solving an embedding problem given a Scholz extension, in the second step we remove ramification in the solution outside the set of exceptional  primes, and in the 
third step 
we force the solution to be Scholz at the cost of one extra ramifying prime. Finally in Section $8$, for $G$ nilpotent this prime is chosen to be the same for all primes $l$ dividing the order of $G$.

We take another approach to the problem of realization of Galois groups with minimal ramification in Section $9$. Take $K = \Q$ or an imaginary quadratic field with $\zeta_l \notin K$. We consider a family of  $l$-extensions of $K$  obtained from central extensions by the Schur multiplicator and observe that a result of Fr\"ohlich \cite{Fro} for $K=\Q$ 
(extended to imaginary quadratic by Watt \cite{Watt}) on realizing the Schur multiplicator confirms Boston's 
conjecture 1.2 of \cite{BM} for groups corresponding to 
this family. 
This conjecture states that for any nontrivial finite group $G$, there exists an extension of $\Q$ with Galois group $G$ and exactly $\max(1, d(G^{ab}))$ ramified primes, and moreover no extension of $\Q$ with Galois group $G$ can be ramified at fewer than $\max(1, d(G^{ab}))$ primes (counting the infinite prime). 
See Kisilevsky-Sonn \cite{kisson} for results on minimally ramified realization of semiabelian groups.

\section{Embedding Problem}
Fix an algebraic closure $\bar K$ of a number field $K$ and let 
 $G_K =
Gal(\bar{K}/K)$ denote the absolute Galois group of $K$. 
An {\emph{embedding problem}} $(G_K, \rho, \alpha)$ for $G_K$ (see e.g. \cite{NSW}) is a
diagram with an exact sequence of finite groups and epimorphism $\rho$. 

\begin{equation}
\cep {C} {G} {\bar G} {\phi} {\rho} {\alpha} \label{cep}\end{equation}

A solution $\phi$ of the embedding problem is a homomorphism $\phi: G_K \rightarrow G$ such that $\alpha \circ \phi = \rho$; a solution is {\emph{proper}} if $\phi$ is surjective. If $G$, $\bar G$ are $l$-groups with the same number of generators, it is easily seen that every solution is proper. 
When the kernel group $C$ is contained in the
center of $G$, the embedding problem ({\ref{cep}}) is called a
{\emph{central embedding problem}}. 
Every nilpotent group can be realized as a Galois group
by solving a sequence of central embedding problems.
For every prime $\p$ of $K$, fix a prime of $\bar K$ above $\p$ and let $D_\p$ (resp. $I_\p$) denote its decomposition (resp.  inertia) subgroup in $G_K$. 

Let 

\begin{equation}
\pcep {C} {G_\p} {\bar G_\p} {\phi_\p} {\rho_\p} {\alpha_\p}  \label{lcep} \end{equation}

denote the corresponding local embedding problem, where $\bar G_\p = \rho(D_\p)$, $G_\p = \alpha^{-1}(\bar{G_\p})$, and $\alpha_\p$, $\rho_\p$ are restrictions of $\alpha$, $\rho$.

In this paragraph we assume in \ref{cep}  that $G$ is an $l$-group and the kernel
$C$ has prime order; let $S_0$ be a finite set  of primes of $K$ containing the
infinite primes, prime divisors of $l$, and prime divisors of a set of
ideals representing the ideal classes of $K$. More generally when $G$ is a
nilpotent group in Section 8, the set $S_0$ contains in addition the
divisors of the order of $G$. It is known, cf. \cite{GJ}, that a solution to a global embedding problem \ref{cep} exists if and only for every prime $\p$ of $K$ there exists a solution to the local embedding problem \ref{lcep}.  The local embedding problem is solvable if $\rho(I_\p)=1$, since $D_\p/I_\p \cong \hat \Z$ is a free group; the Scholz condition ensures solvability at the ramified primes. Let $Ram(\rho) = \{ \p {\rm  \hs of \hs } K \mid \rho(I_\p) \not = 1\}$. 

\begin{defn}[$l^N$-Scholz homomorphism, cf. 3.2 of \cite{GJ}]

Given a number field $K$, an $l$-group $G$ and a positive integer $N$ 
such that $l^N$ is divisible by the exponent of $G$. Denote by $T$ a set of $l^N$-exceptional primes as defined in Section 4. 

An epimorphism $\phi:G_K \rightarrow G$ is $l^N$--Scholz if
\begin{itemize}
\item For $\p \in Ram(\phi) \cup T$, $\phi(D_\p) = \phi(I_\p)$.
\item For $\p \in Ram(\phi)$, the absolute norm $ N(\p) \equiv 1 \pmod {l^N}$. 
\item For $\p \in S_0$, $\phi(D_\p)=1$.
\label{mScholzHom}
\end{itemize}
The last condition is an example of local data of \cite{GJ}. We will also say the extension $L/K$ is $l^N$-Scholz, where $L$ is the subfield of $\bar K$ fixed by $ker(\phi)$.

\end{defn}

The definition of $l^N$-Scholz does not depend on the choice of prime of $\bar K$ above each $\p$. Clearly  if a homomorphism $\phi$ is $l^N$-Scholz, then it is $l^k$-Scholz for all integers $k \leq N$. 

\section{Existence of Solutions}

\begin{thm}[Existence Theorem]
Let $(G_K, \rho, \alpha)$ be a central embedding problem, $\bar G = \rho(G_K)$ is an $l$-group, $C = ker(\alpha)$ cyclic of order $l^e$. Suppose $\rho$ is $l^N$-Scholz (exponent of $G$ divides $l^N$) 
and $\zeta_l \notin K$. Then the embedding problem

\begin{equation}
\cep {C} {G} {\bar G} {\psi_0} {\rho} {\alpha} \label{ETcep} \end{equation}

has a solution. 
\label{ExistLemma}\end{thm}
 
\proof
If $G$ is a split extension of $\bar G$ , we may apply Proposition \ref{SplitCase}, so assume the extension is Frattini, i.e., $C$ is contained in the Frattini subgroup of $G$.  
We may break \ref{ETcep}  into a sequence of $e$ embedding problems each with kernel group of  order $l$, which we may solve by Proposition 7.3 of \cite{GJ} at the cost of one ramified prime at each step. We obtain an $l^N$-Scholz solution $\psi_0$ to \ref{ETcep} such that

\[Ram(\psi_0) \cup T = Ram(\rho) \cup T \cup \{e {\rm { \hs primes \hs of \hs }} K\}\]
\done

In sections 5-7 we will show that the embedding problem \ref{ETcep} has an $l^N$-Scholz solution 
at the cost of only one additional ramified prime (assuming $K$ has
no ideal classes of order $l^2$ if $|C| > l$).


\section{Exceptional Set of Primes}

The key Lemma \ref{FieldIsoLemma} was originally proved by the first author in a different way in her thesis
\cite{MarkinThesis}. 
The lemma below generalizes results of Gras in \cite{GrasCFT} Ch. II, Theorem 6.3.2 and Lemma 4.1, p. 361 in \cite{RubinECCM}. 
\begin{lem} \label{KL}

Let $L/K$ be a Galois $l$-extension,  $\tilde K = K(\mu_m)$, $\tilde L = L(\mu_m)$,  where $m$ is a power of $l$. If $\zeta_l  \notin K$, then the canonical map
\[\Km/\Km^m \rightarrow \tilde L^\times/\tilde {L^\times}^m \]
is injective. \label{KumInj}
\end{lem}

\proof
From Kummer theory, we have 
$H^1(Gal(\bar K/K),\mu_m) \cong \Km/ {\Km}^m$ and 
$\\H^1(Gal( \bar K/\tilde{L}),\mu_m) \cong \tilde L^\times/\tilde{L^\times}^m,$
 where $\bar K$ denotes an algebraic closure of $K$. 
 The extensions $K \subseteq \tilde L \subseteq \bar K$ give the following exact sequence of cohomology groups via the restriction-inflation maps
\[ \les {H^1(Gal(\tilde L /K), \mu_m^\gamma)} {H^1(Gal(\bar K/K), \mu_m)}  {H^1(Gal(\bar K/\tilde L),\mu_m)},\]
 where $\gamma = Gal(\bar K/\tilde L)$.  It suffices to prove $H^1(Gal(\tilde L/K),\mu_m^\gamma)=0$; note $\mu_m^\gamma = \mu_m$. By a second application of the restriction-inflation sequence, now to the extensions $K \subseteq L \subseteq \tilde L$, we have the exact sequence

\[ \les {H^1(\Gamma/\Delta,\mu_m^\Delta)} {H^1(\Gamma, \mu_m)}{H^1(\Delta, \mu_m)},\]

where $ \Gamma = Gal(\tilde L/K)$, $ \Delta = Gal(\tilde L / L)$. 
The cohomology group ${H^1(\Gamma/\Delta,\mu_m^\Delta)} = 0$ since $\mu_m^\Delta = \mu_m \cap L = \{1\}$ ( $\zeta_l \notin K$ and $L/K$ is an $l$-extension). 
Since $\Delta$ is cyclic, by Herbrand theory, the orders of the Tate cohomology groups $H^i(\Delta, \mu_m)$ are equal for $i=0,1$. 
But $H^0(\Delta, \mu_m) = \mu_m^\Delta/Norms = 0$. This completes the proof. 
\done

Let $K_S$ be the group of $S$-units of $K$, where $S$ contains the infinite primes of $K$. By Dirichlet's unit theorem, the $\Z$-rank of  $K_S$ is $$u:=rk_{\Z}(K_S)=|S|-1.$$

\begin{lem}

Assume $\zeta_l \notin K$. 
With the notation of Lemma  \ref{KumInj} let $M $ be an abelian extension of $L$ containing $\tilde L$. 
There are isomorphisms 
\[Gal(\tilde K(\rmKS)/\tilde K) \stackrel{f_1}\cong Gal(\tilde L(\rmKS)/\tilde L) \stackrel{f_2}\cong Gal(M(\rmKS)/M) \cong (\Z/m\Z)^u.\] \label{isomexts}
\label{FieldIsoLemma}
\end{lem}

\proof
 
Apply Lemma \ref{KumInj} restricted to the image of $K_S$ in $\tilde L^\times$ to conclude that $f_1$ is an isomorphism.
Next we show $f_2$ is an isomorphism. Let $F = \tilde L(\rmKS) \cap M$. We show $F = \tilde L$, so that $f_2$ would be an isomorphism.  Since $F \subset M$, the extension $F/L$ is abelian. And $\tilde L \subset F \subset L(\rmKS)$.
If $F$ is not $\tilde L$, then $F$ contains a cyclic extension $F_0/\tilde L$, $[F_0: \tilde L]=l$. From Kummer theory, $F_0 = \tilde L(\sqrt[l]{b})$, $b \in K_S$. But $Gal(\tilde L(\sqrt[l]{b})/L)$ is not abelian, thus $F = \tilde L$. 
\done

The corollary below will be used in section \ref{NilpotentSection}. We include it here for convenience. 

\begin{cor}
\label{liftinglemma} Let $K$ be a number field,  $S$ a finite set of primes of $K$ and let $a>1$ be an integer. For each $l \mid a$ let $L_l /K$ be a Galois
$l$-extension. Suppose that $\zeta_{l} \not\in K$ for each $l \mid a$. Set
$M_l = L_l(\zeta_{l^{N}}, \zeta_a)$ and $ M = \prod_l M_l$.  Then we have a series of isomorphisms:

$$Gal(K(\rKSl)/K(\zeta_{l^{N}})) \cong Gal(L_l(\rKSl)/L_l(\zeta_{{l}^{N}}))$$ 
$$ \cong Gal(M_l(\rKSl)/M_l)\cong Gal(M(\rKSl)/M)  $$

\end{cor}

The diagram below contains the fields involved in these isomorphisms. 
\[
           \xymatrix{K( \rKSl)\ar[r] &L_l(\rKSl) \ar[r] &M_l(\rKSl)\ar[r]   &M(\rKSl)           \\
           K(\zeta_{l^{N}})\ar[r]\ar[u] &L_l(\zeta_{{l}^{N}})\ar[r]\ar[u] &M_l\ar[r]\ar[u] &M\ar[u]}
               \]

\proof  

The first two isomorphisms follow from Lemma \ref{FieldIsoLemma}. To show the rightmost isomorphism 
note that $M_l(\rKSl)/M_l$ is an $l$-extension, while $l \nmid [M : M_l]$. \done

\begin{lem}

For each $l | a$, assume that $\zeta_l \notin K$. Let $\N_l$ denote the field $ L_l(\rKS)$ and let $\sigma_l \in Gal( \N_l/L_l(\zln))$. Define $\N = \prod_{l | a}\N_l$. Then there exists ${\sigma} \in Gal(\N/K(\mu_a))$ such that
${\sigma}|_{\N_l} = \sigma_l$ for all $l | a$.
\label{corlift}
\end{lem}

\proof
By Lemma \ref{liftinglemma}, each $\sigma_l$ extends to an element, say $\hat \sigma_l$, 
of $Gal(\N_lM_l/M_l)$. The latter group is a subgroup of the $l$-group
$Gal(\N_lM_l/K(\mu_a))$. Now observe that $Gal(\N/K(\mu_a)) \cong \prod_{l | a}
Gal(\N_lM_l/K(\mu_a))$. Therefore we may define  $\sigma \in Gal(\N/K(\mu_a))$ as $\sigma = \prod_{l \mid a}{\hat \sigma_l}$.

\done

For an abelian group $A$ and a prime number $l$, let $A_l = \{a \in A \mid a^l=1\}$. 
We define the subgroup $V(l)$ (denoted by $V$ when the prime $l$ is implicit) of $\Km$ $$V = \{ a \in K^\times \mid (a) = \mathfrak a ^l \hs \text{for  a fractional ideal $\mathfrak a$ of $K$} \}.$$ We have the following split exact sequence, e.g. pg. 109 of \cite{KochGalpGerman},
$$ \es {E/E^l} {V/{K^{\times l}}} {Cl(K)_l}, $$

where $E$ denotes the group of units of $K$ and the right hand map sends $\\a \mod K^{\times l}$  to the ideal class of $\mathfrak a$, where $(a) = \mathfrak a ^l$. Similarly 

$$ \es {E/E^{l^N}} {EV^{l^{N-1}}/{K^{\times {l^N}}}} {Cl(K)_l}.$$

Let $w_1, \ldots, w_s$ be a $\Z$-basis of $E$ mod torsion. As in \cite{GJ}, choose ideles $\alpha_1, \ldots, \alpha_r \in J$ whose images are an $\F_l$-basis of the $l$-torsion subgroup $(J/{\Km U})_l$ of the ideal class group of $K$. Then for $j=1, \ldots, r$
$$\alpha_j^l=a_j^{-1}\epsilon_j, \hs \ \ a_j \in \Km, \epsilon_j = (\epsilon_{j,v}) \in U, \epsilon_{j,v} \in U_v.$$

For all $j$ and all primes $v$ of $K$, $a_j$ and $\epsilon_{j,v}$ have the same image in $U_v/{U_v^l}$. Taken $\mod {\Km}^l$, the set 
$\{w_1, \ldots , w_s, a_1, \ldots, a_r\}$ is a basis of $V/{{\Km}^l}$.

We define a governing field $\Omega_l$ (compare Chapter 5 of \cite{GrasCFT} or \cite{GJ} for $N=1$)
\begin{equation}
\label{omega}
\Omega_l = K(\mu_{l^N}, \sqrt[l^N]{EV^{l^{N-1}}}) = 
K(\mu_{l^N}, \sqrt[l^N]{E}, \sqrt[l]{V}) = \end{equation}
\[K(\mu_{l^N}, \sqrt[l^N]{w_i}, \sqrt[l]a_j: {{1 \leq i \leq s, 1 \leq j \leq r}} ). 
\]

It follows from lemma \ref{isomexts} that the Kummer extension satisfies $Gal(\Omega_l/K(\mu_{l^N})) \cong (\Z/{l^N \Z})^s \oplus   (\Z/{l \Z})^r.$ Of course, if $K = \Q$, we have $r = s = 0$.

Define subfields of $\Omega_l$
\[N_i = K(\mu_{l^N}, \sqrt[l^N]{w_{k}} : { 1 \leq k \leq s, k \not = i ; \hs \sqrt[l]{a_k} : \hs  1 \leq k \leq r} ), \hs  1 \leq i \leq s \]

\[N_j' = K(\mu_{l^N}, \sqrt[l^N]{E}; \hs \sqrt[l]a_{k} : 1 \leq k \leq r, k \not = j ), \hs 1 \leq j \leq r. \]

Then $Gal(\Omega_l/N_i)$ is cyclic of order $l^N$ while $Gal(\Omega_l/{N_j'})$ has order $l$.

\begin{defn}[cf. (5.5) of \cite{GJ} for $N=1$] 
\label{ExcPrimes}

A set $T_l = \{\p_1, \ldots, \p_s, \q_1, \ldots \q_r\}$ of prime ideals of $K$ such that $T_l  \cap S_0 = \emptyset$ is {\emph{$l^N$-exceptional}} if

$$Gal(\Omega_l/N_i) = D_{\p_i}(\Omega_l/K), \hs 1\leq i \leq s \hs {\rm{and}} $$

$$Gal(\Omega_l/N_j') = D_{\q_j}(\Omega_l/K), \hs 1\leq j \leq r.$$ 

Note that this property is independent of the primes above $\p_i$  (resp. $\q_j$) since $N_i$ (resp. $N_j'$) is a normal extension of $K$. 

\end{defn}
For a prime ideal $\p$ of $K$ unramified in a Galois extension $F/K$, $\Frob {F/K}{\p}$ denotes the conjugacy class in $Gal(F/K)$ consisting of the Frobenius elements of all prime ideals of $F$ above $\p$.

Choose $$\sigma_i(l) \in \Frob {\Omega_l/K} {\p_i(l)}, \hs 1\leq i \leq s \hs {\rm and} $$
  $$\tau_j(l) \in \Frob {\Omega_l/K} {\q_j(l)}, 1\leq j \leq r_l;$$ here we make the dependence on $l$ explicit. Note that  $\{\sigma_i(l), \tau_j(l): 1\leq i \leq s, 1\leq j \leq r_l\}$ are a minimal generating set of the abelian group $Gal(\Omega_l/K(\mu_{l^N}))$. Further if $a$ is the product of the primes dividing $|G|$, the latter group is isomorphic to $Gal(\Omega_l(\mu_{a^N})/K({\mu_{a^N}}))$ by Lemma \ref{FieldIsoLemma}.

By the Chebotarev density theorem, there exists an $l^N$-exceptional set of primes disjoint from 
any given set of primes of $K$ of density $0$. Note that since $v$ splits completely in $K(\mu_{l^N})/K$ for all $v \in T_l$, we have $\zeta_{l^N} \in K_v$ for all $v \in T_l$.

It follows from Kummer theory for primes $\p_i, \q_j \in T_l$ that

\begin{itemize}

\item $w_i$ not an $l$-th power in $U_{\p_i}$
\item  $w_i \in U_v^{l^N}  \hs \forall v \in T_l,  \hs v \not = \p_i$.

\vspace{.5 cm}

\item $a_j$ not an $l$-th power in $U_{\q_j}$
\item $a_j \in U_v^{l}  \hs \forall v \in T_l, \hs v \not = \q_j$.

\end{itemize}

Note that if $T_l$ is $l^N$-exceptional, then $T_l$ is $l^k$-exceptional for all $1 \leq k \leq N$. We will therefore fix a set $T_l$ of $l^N$-exceptional primes, where $l^N$ is divisible by the exponent of the $l$-group $G$. From now on until section 8 we will let $T$ denote $T_l$, as the prime $l$ is implicit.

\section{Split Case}

We begin with a lemma which generalizes  Lemma 4.2 of \cite{GJ}. 
If $K = \Q$, and $b$ is an integer greater than one, the lemma follows at once from the fact there are infinitely many primes $q \equiv 1 \pmod b$, and we take subfield $M$ of $\Q(\mu_q)$ of degree $b$.

\begin{lem}
Given integer $b > 1$ and number field $K$, there exist infinitely many prime ideals $\q$ of $K$ and cyclic extensions $M = M(\q)$ of $K$ of degree $b$ such that $\q$ is the unique ramified prime of $M/K$, $\q$ is totally ramified, and $\q$ does not divide $b$. 
\label{qram}
\end{lem}
\proof

Let $S$ be a finite set of primes of $K$ containing $S_0$ and prime divisors of $b$ and let $\Omega = K(\sqrt[b]K_S)$. By Chebotarev's theorem there exist infinitely many primes $\q$ of $K$, $\q \notin S$, such that $\q$ splits completely in $\Omega/K$. For such $\q$, $\Omega$ is contained in the completion $K_\q$ and so $K_S \subset (K_\q^\times)^b$.

Define $J_S = \prod_{v \in S}K_v^\times \times \prod_{v \notin S}U_v \subset J$. By class field theory, cyclic extensions of $K$ are given by idele class characters. Since $J/K^\times \cong J_S/K_S$, we want to define an epimorphism $\chi:J_S/K_S \rightarrow \mu_b$ with $\chi(K_S) = \{1\}$. The group $U_\q/U_\q^b$ is cyclic of order $b$, so there is an epimorphism $\chi_\q:U_\q \rightarrow \mu_b$ with kernel $U_\q^b$. For $\alpha = (\alpha_v) \in J_S$, define $\chi(\alpha)=\chi_\q(\alpha_\q)$. Note $\chi(K_S) = \{1\}$ and $\chi(K_v^\times) = \{1\}$, $v \in S$. By class field theory, $\chi$ corresponds to a cyclic, degree $b$ extension $M(\q)/K$ in which $\q$ is totally and tamely ramified and the other primes of $K$ are unramified. 
\done

\begin{thm}

Let $A$ be a finite abelian group with $d$ generators. There exist infinitely many Galois extensions $N/K$ such that $Gal(N/K) \cong A$ and exactly $d$ primes of $K$ ramify in $N$. Such $N$ is its own genus field relative to $K$. 

\label{AbG}
\end{thm} 

\proof

Write $A$ as a direct product of $d$ cyclic groups and apply Lemma \ref{qram} to each factor. The resulting extensions $M(\q_i)$, $1 \leq i \leq d$ are linearly disjoint over $K$ by ramification considerations. Take $N$ to be the composite of the fields $M(\q_i)$. 
Note that these $\q_i's$ are not to be confused with the ones defined in Definition \ref{ExcPrimes}. \done

\begin{prop}[Split Case]
Let  $G$ be an $l$-group of exponent dividing $l^N$. Suppose the homomorphism $\rho: G_K \rightarrow \bar G$ is $l^N$-Scholz and the central exact sequence is split
\[\es C G {\bar G},\] where the kernel $C$ of $\alpha: G \rightarrow \bar G$ is cyclic.
There is an $l^N$-Scholz solution $\phi$ to the embedding problem $(G_K, \rho, \alpha)$ and a prime $\q $ not in $S = Ram(\rho) \cup S_0 \cup T$ such that $Ram(\phi) = Ram(\rho) \cup \{\q \}$. \label{SplitCase}

\end{prop}

\proof

We apply the argument in Lemma \ref{qram} with $b = |C|$, $\Omega = L(\mu_{l^N}, \sqrt[b]{K_S})$, where $L$ is the subfield of $\bar K$ fixed by $ker(\rho)$,  to obtain $\q$ and an idele class character $\chi$ of order $b$; $\q$ splits completely in $\Omega/K$. By the Reciprocity law $\chi$ corresponds to an epimorphism $\eta: G_K \rightarrow C$. Then $\phi = (\rho, \eta): G_K \rightarrow \bar G \times C$, $\sigma \mapsto (\rho(\sigma), \eta(\sigma))$, is a proper solution to the embedding problem. It remains to check that $\phi$ is $l^N$-Scholz, given that $\rho$ is $l^N$-Scholz. 

If $v \in S_0$, $\phi(D_v) = 1$ since $\rho(D_v)=1$ (given) and $\eta(D_v)=1$ for $v \in S$. If $v \in T$, $\phi(D_v)= \phi(I_v)$ since $\rho(D_v)=\rho(I_v)$ (given) and $\eta(D_v) = 1$ for $v\in S$. 

Suppose $v \in Ram(\phi) = Ram(\rho) \cup \{\q\}$. 

If $v = \q$, $\q$ splits completely in $K(\mu_{l^N})/K$, hence $N(\q) \equiv 1 \pmod {l^N}$. Since $\q$ splits completely in $L/K$, $\rho(D_\q)=1$. As $\eta(I_\q)=C$ for $\eta:G_K \rightarrow C$, we have $\eta(D_\q)=\eta(I_\q)$. Thus $\phi(D_\q)=\phi(I_\q)$. 

If $v \in Ram(\rho)$, then $N(v) \equiv 1 \pmod {l^N}$ and $\rho(D_v) = \rho(I_v)$ (given). But $\eta(D_v) = 1$ since $v \in Ram(\rho) \subset S$. Thus $\phi(D_v)= \phi(I_v)$ for $v \in Ram(\phi)$. 

We conclude $\phi = (\rho, \eta)$ is an $l^N$-Scholz solution with one additional ramified prime. 
\done

\section{Removing Ramification}
\begin{lem}
Let $K$ be a number field not containing $\zeta_l$, $N \geq e \geq 1$. Given a finite set $S$ of primes disjoint from an $l^N$-exceptional set $T$, characters $\chi_v: U_v \rightarrow \mu_{l^e}$, for $v \in S$, at least one of which is onto. Assume $K$ has no ideal classes of order $l^2$ when $e > 1$. There exists an idele class character 
$\chi: J/K^\times \rightarrow \mu_{l^e}$ such that 
$\chi|_{U_v} = \chi_v$ for all $v \in S$ and $\chi|_{U_v} = 1$ for all $v \notin S \cup T$.

\label{remramlem}

\end{lem}

\proof

It suffices to prove the result when $S=\{v_0\}$ and then take the product of the resulting characters. 
Let $I=T \cup \{v_0\}$. 

{\bf {Step 1: Defining $f$ on $U\Km/\Km$.\\}} 
We define an epimorphism $f: U \rightarrow \mu_{l^e}$ of the form

\[f= \prod_{v \in I}\chi_v ,\]

with $f|_{U_v}=1$ for $v \notin I$. 
The character  $\chi_{v_0}$ is given and the characters $\chi_v, \hs v \in T,$ are to be defined suitably. Each character $\chi_v$ is trivial for $v \notin I$. 
 
By the definition of an $l^N$-exceptional set of primes, the image of each unit $w_i$ generates  $U_{\p_i}/{U_{\p_i}^{l^e}}$, $\p_i \in T$,  hence we can define $\chi_{\p_i}: U_{\p_i} \rightarrow \mu_{l^e}, \hs 1 \leq i \leq s,$ to satisfy $$\chi_{\p_i}(w_i)\chi_{v_0}(w_i) =1.$$

Similarly  $\epsilon_{j,\q_j}$ generates $U_{\q_j}/U_{\q_j}^{l}$ (hence also modulo $U_{\q_j}^{l^e}$) and we can define 
$\chi_{\q_j}: U_{\q_j} \rightarrow \mu_{l^e}, \hs 1 \leq j \leq r,$  to satisfy $$\chi_{\q_j}(\epsilon_{j,\q_j})\chi_{v_0}(\epsilon_{j,{v_0}}) =1.$$

Next we establish the "off-diagonal" vanishing of $\prod_{v \in I}\chi_v$. Recall that $\epsilon_{j,v} \in U_v^{l}$ for $\q_j \not = v \in T$ for each $j$, and  $w_i  \in  U_v^{l^e}$ for $\p_i \not = v \in T$ for each $i$. Thus we have

$$\prod_{v \in I}\chi_v(w_i) = \chi_{\p_i}(w_i) \chi_{v_0}(w_i)\prod_{\p_i \not = v \in T}\chi_v(w_i) = 1,$$

$$\prod_{v \in I}\chi_v(\epsilon_{j,v}^{l^{e-1}}) = \chi_{\q_j}(\epsilon_{j,\q_j}^{l^{e-1}}) \chi_{v_0}(\epsilon_{j,v_0}^{l^{e-1}})\prod_{\q_j \not = v \in T}\chi_v(\epsilon_{j,v}^{l^{e-1}}) = 1.$$

It follows that $\prod_{v \in I}\chi_v$ is trivial on the image of 
$E \oplus (\oplus_{j=1}^r \langle \epsilon_j \rangle ) $ in $\prod_{v \in I}U_v/{U_v^{l^e}}$. Letting $\Delta: K^\times \rightarrow J$ be the diagonal embedding, we have in particular $f(\Delta(E))=1$, so $f$ is defined on $U/\Delta(E)$, which we write as $U/E \cong U\Km/\Km$. 

Note if $l$ does not divide the class number of $K$, then $f$ already provides
the desired idele class character since the $l$-part of the ideal class group $J/\Km U$ will be
trivial. Otherwise we must extend $f$ from $\Km U/\Km$ to $J/\Km$. 

{\bf{Step 2: Character of order $l$. \\}}
Define 
$f_1: U \rightarrow \mu_l$ by $f_1 = f^{l^{e-1}}$. By the techniques  of the proof of Lemma 6.1 of \cite{GJ}, $f_1$ extends to an idele class character $\chi_1$ of order $l$ with $\chi_1|_{U_{v}}=\chi_{v}^{l^{e-1}}$, for $v \in I$ and $\chi_1|_{U_v}=1$ if $v \notin I$. This follows from the trivial fact that an $l^e$-exceptional set $T$ is $l$-exceptional.

We have ${\frac{\Km U}{\Km ker(f_1)}} \cap {\frac{\Km ker(\chi_1)}{\Km ker(f_1)}} \equiv 1$. 
Also we have 
$\\|J/{\Km ker(f_1)}| =|J/\Km U| \cdot |\Km U/\Km ker(f_1)|= h \cdot l$, where $h$ is the class  number of $K$, which we may assume is
a power of $l$.
Thus
$|\Km ker(\chi_1)/ \Km ker(f_1)|= \frac{|J/{\Km ker(f_1)}|}{|{J/{\Km ker(\chi_1)}}|} =\frac{h \cdot l}{l} = |J/{\Km U}|$. This implies that the following exact sequence 

$$\es {\frac{\Km U}{\Km ker(f_1)}} {\frac{J}{\Km ker(f_1)}} {J/{\Km U}}$$

splits, with  $ {\frac{\Km ker(\chi_1)}{\Km ker(f_1)}}$ mapping isomorphically onto ${J/{\Km U}}$. The image
 ${J/{\Km U}}$ has exponent $l$ by assumption and the kernel is cyclic of order $l$. Hence ${\frac{J}{\Km ker(f_1)}}$ has exponent $l$. 

{\bf{Step 3: Extending to a character of order $l^e$.\\}}
First we prove the following claim about finite abelian $l$-groups. 
\begin{claim}\label{dirsum}
Let $\Gamma$ be a finite abelian $l$-group and let  $\gamma \subseteq \Gamma$ be a cyclic subgroup of order $l^e$. If $\Gamma/{\gamma^l}$ has exponent $l$, then $\gamma$ is a direct summand of $\Gamma$.

{\emph{Proof:\\}}
The exponent of $\Gamma$ is $l^e$, since for any element $g \in \Gamma$ we have $g^l \in \gamma^l$ and hence $g^{l^e} =1$. Therefore $\gamma$ is a subgroup generated by an element of maximal order, and hence is a direct summand, as desired. 
\end{claim}

We have the following diagram with exact rows and columns

\begin{xymatrix}{
&1 \ar[d] &1 \ar[d] \\
& \frac{(\Km U)^l \Km ker(f)}{\Km ker(f)} \ar[r]^{=} \ar[d] & \frac{(\Km U)^l \Km ker(f)}{\Km ker(f)}  \ar[d] \\
1 \ar[r] & \frac{\Km U}{\Km ker(f)} \ar[r]\ar[d] &  {J/{\Km ker(f)}} \ar[r]\ar[d]^{} & {J/{\Km U}} 
 \ar[r]\ar[d]^{\cong} &1\\
 1 \ar[r] & \frac{\Km U}{\Km ker(f_1)} \ar[r]\ar[d] & {J/{\Km ker(f_1)}}  \ar[r]\ar[d] & {J/{\Km U}}  \ar[r] &1\\
&1  &1
}
\end{xymatrix}
We apply Claim \ref{dirsum} with $\gamma =\frac{\Km U}{\Km ker(f)}$ and $\Gamma = \frac{J}{\Km ker(f)}$. It follows from the diagram above that $\Gamma/{\gamma^l} \cong \frac{J}{\Km ker(f_1)}$, which by assumption has exponent $l$. Claim \ref{dirsum} implies that $\frac{\Km U}{\Km ker(f)}$ is a direct summand of   
$\frac{J}{\Km ker(f)}$. Thus $f$ extends to a character $$\chi: J/\Km \rightarrow \mu_{l^e}$$ by defining $\chi$ equal to $f$ on $U$ and $\chi$ trivial on a complement of $\Km U/{\Km ker(f)}$.
\done

\begin{thm}[Removing Ramification]
\label{remram}

Suppose $K$ has no ideal classes of order $l^2$ and does not contain $\zeta_l$. If the Frattini embedding problem 
 $(G_K, \rho, \alpha )$ has a solution $\psi_0$, then it has a solution $\psi:G_K \rightarrow G$ with $Ram(\psi) \subset Ram(\rho) \cup T$.

\end{thm}
\proof

The proof is similar to Lemma 6.2 of \cite{GJ} except that we twist $\psi_0$ by a character of order $l^e$. Let $S = Ram(\psi_0) \setminus \{Ram(\rho) \cup T\}$, so if $v \in S$, then $\psi_0(I_v) \subseteq C$. Set $l^e = \max \{ |\psi_0(I_v)|: v \in S \}$.

For $v \in S$ we define $\chi_v:= \psi_0|_{I_v}$ viewed as $\chi_v: U_v \rightarrow \mu_{l^e}$ by reciprocity. By  \ref{remramlem} there exists an idele class character $\chi$ of order $l^e$ with certain local properties. We identify $\chi$ with $\eta: G_K \rightarrow C$ via reciprocity and set $\psi = \psi_0 \eta^{-1}$. Since the embedding problem $(G_K, \rho, \alpha )$ is Frattini, $\psi$ is surjective. 

\done

\begin{rem}
Note that in case $e=1$ the hypothesis on the order of ideal classes in the theorem above 
can be dropped. 
\end{rem}

\section{Finding an $m$-Scholz solution}

We generalize Lemma 7.1 of  \cite{GJ} to prime powers. 

\begin{lem}
Given integers $N \geq e \geq 1$, Galois $l$-extension $L/K$, characters $\chi_v : K_v^\times \rightarrow \mu_{l^e}$ for all
$v$ in a finite set $S \supseteq S_0$. 
Assume that $K$ does not contain $\zeta_l$. 
There exists a prime ideal $\q$ of $K$ outside $S$ and a character $\chi: J_K/{K^\times} \rightarrow \mu_{l^e}$ such that conditions (1)-(4) hold:
\begin{itemize}
\item $\q$ splits completely in $L(\mu_{l^N})/K$
\item
$\chi|_{\Kvs}=\chi_v$ for all $v \in S$. 
\item

$\chi(U_\q) = \mu_{l^e}$.

\item $\chi(U_v) = 1 $ for all $ v \notin S \cup \{\q\}$. 
\end{itemize}

\label{mScholzSoln}
\end{lem}

\proof

Since $S_0$ is chosen large enough, we have $J_S/K_S \cong J/{\Km}$. It therefore suffices to define a character $g: J_S \rightarrow \mu_{l^e}$ such that for all $(\alpha_v) \in J_S$ 
\[g((\alpha_v))  = \chi_\q(\alpha_\q) \times \prod_{v \in S}\chi_v(\alpha_v)\]

for some prime $\q$ and some epimorphism $\chi_\q: U_\q \rightarrow \mu_{l^e}$ chosen so that $\q$ splits completely in $L(\mu_{l^N})/K$ and $g(K_S) =\{1\}$.

We define a character $h:K_S \rightarrow \mu_{l^e}$ as the composition

\[K_S \stackrel{j}{\rightarrow}  J_S  \rightarrow  \mu_{l^e} \]
 where the left map $j$ is the embedding of $K_S$ in $\prod_{v \in S}K^{\times}_v$ and the right map is ${\prod_{v \in S}\chi_v}$. Thus for $x \in K_S$, $g(x) = h(x)\chi_\q(x)$, so $\chi_\q$ must be chosen to make $g(x)=1$ for all $x \in K_S$. 

{\bf \emph{Case $h(K_S)=\{1\}$}}. If $\q$ satisfies  $K_S \subset U_\q^{l^e}$, then for any character $\chi_\q: U_\q \rightarrow \mu_{l^e}$, we have $\chi_\q(K_S)=\{1\}$.
By Chebotarev's theorem, there exists a prime ideal $\q \notin S$ of $K$ which splits completely in $\Omega := L(\mu_{l^N}, \sqrt[l^e]{K_S})$. Note that $\q$ splitting completely in $K(\mu_{l^N})/K$ implies that absolute norm $N^K_\Q(\q) \equiv 1 \pmod {l^N}$. Then $K_S \subseteq U_\q^{l^e}$ by Kummer theory.

{\bf \emph{ Case $h(K_S) \neq \{1\}$}}. The image $h(K_S)$ is cyclic of order $l^k$, $1 \leq k \leq e$. Thus there exists  $x_1 \in K_S$ with $h(x_1)$ of order $l^k$.  $K_S/{K_S^{l^k}}$ may be generated by $ \{ x_1, x_2, \ldots, x_u \} $, with $h(x_i) =1 , i > 1$. By Burnside's basis theorem $\{x_1, \ldots , x_u\}$ also generate $K_S/K_S^{l^e}$. We want to pick a prime $\q \nmid l$, $\q \notin S$ such that

\begin{itemize}
\item $\q$ splits completely in $L(\mu_{l^N})/K$.
\item $x_1 \in U_\q^{l^{e-k}} \setminus U_\q^{l^{e-k+1}}$.
\item $x_i \in U_\q^{l^e}$ if $i > 1 $.
\end{itemize}

To that end let 
\[\Omega_k = L(\mu_{l^N}, \sqrt[l^{e-k}]{x_1}, \sqrt[l^e]{x_i} : {i > 1}).\]

The field $\Omega_k$ is a normal extension of $K$. By Lemma \ref{FieldIsoLemma}, $Gal(\Omega/L(\mu_{l^N})) \cong (\Z/{l^e}\Z)^u$  and $Gal(\Omega/\Omega_k)$ is cyclic of order $l^k$. By Chebotarev's theorem we may choose $\q \notin S$ such that $Frob(\q, \Omega/K)$ generates $Gal(\Omega/\Omega_k)$, in particular $\q$ splits completely in $\Omega_k/K$. This guarantees that the above three conditions on $\q$ are satisfied. 

Having chosen $\q$, we define $\chi_\q$, a character of order $l^e$. 
Choose $y \in U_\q$ such that $y^{l^{e-k}} = x_1 \in U_\q$. We want $\chi_\q(y)$ of order $l^e$, then $\chi_\q(x_1)$ has order $l^k$. If $\beta = h(x_1)$ is an element of $\mu_{l^e}$ of order $l^k$, then $\beta = \alpha^{l^{e-k}}$, where $\alpha$ is a generator of $\mu_{l^e}$. Set $\chi_\q(y) = \alpha^{-1}$. Then $\chi_\q(x_1)= \beta^{-1}$.

So we have chosen $\chi_\q$ so that $\chi_\q(x_1)h(x_1)=1$. Thus $g(K_S)=1$ and we have proved the lemma for prime power order characters.  \done

\begin{prop}
Suppose that the central embedding problem $(G_K,  \rho, \alpha)$, $G$ an $l$-group, is Frattini, $\rho$ is $l^N$-Scholz, and $\zeta_l \notin K$. Assume there
exists a solution $\psi$ with $Ram(\psi)\cup T = Ram(\rho) \cup T$. Then there exists a prime $\q \notin S:=Ram(\psi)\cup S_0\cup T$ and an $l^N$-Scholz solution $\varphi$ such that $Ram(\varphi)=Ram(\psi)\cup\{\q \}$. \label{nScholzSol}

\end{prop}
\proof

{\bf\emph{Step 1.}} Define homomorphisms $\eta_v: D_v \rightarrow C$, $v \in S$. 

$\bullet$ If $v \in S \setminus S_0$, we lift Frobenius at $v$ to $\sigma_v \in D_v$. Since $\rho$ is $l^N$-Scholz and $Ram(\psi) \cup T = Ram(\rho) \cup T$, after adjusting the lift $\sigma_v$ we may assume $\psi(\sigma_v) \in C$ (see pg. 36 of \cite{GJ}). Then let 
$\eta_v$ be the unique homomorphism $D_v \rightarrow C$ satisfying $\eta_v(\sigma_v) = \psi(\sigma_v)$ and $ \eta_v(I_v) = \{1\}$. 

$\bullet$ If $v \in S_0$, $\alpha(\psi(D_v)) = \rho(D_v)=\{1\}$, again since $\rho$ is $l^N$-Scholz. Thus $\psi(D_v) \subset ker(\alpha) = C$. So define $\eta_v = \psi|_{D_v}$. 

We have defined $\eta_v$, $v \in S$; now we apply Lemma \ref{mScholzSoln} to get a map $\eta: G_K \rightarrow C$ and a prime $\q \notin S$ such that $\eta|_{D_v}= \eta_v$, $v \in S$, $\eta(I_\q) = C$, and $\eta$ unramified for $v \notin Ram(\psi) \cup T \cup \{\q\}$. Finally set $\varphi = \eta^{-1}\psi$. Note $\varphi(\sigma_v)=1$, so $\varphi(D_v) = \varphi(I_v)$ if $v \in Ram(\psi) \cup T \setminus S_0$.

{\bf \emph{Step 2.}} We claim $\varphi$ is unramified outside $Ram(\psi)\cup\{\q\}$. In fact if $v \in S\setminus S_0$, we have $\eta(I_v)=\eta_v(I_v)=\{1\}$, so $\varphi(I_v)=\psi(I_v)$. The result follows.

{\bf \emph{Step 3.}} We claim $\varphi$ is $l^N$-Scholz. Since the extension is Frattini, any solution is proper. The check of the three points of Definition \ref{mScholzHom} is similar to pg. 37 of \cite{GJ} except for the proof that $\varphi(D_\q) = \varphi(I_\q)$. For that, note that $\q$ is chosen to split completely in the fixed field of $ker(\psi)$, so $\psi(D_\q)=\{1\}$. Putting this together with $\eta(I_\q)=C$, we conclude that $\varphi(D_\q) = \varphi(I_\q)$.

\done

Putting together Existence Theorem, Proposition \ref{SplitCase}, Proposition \ref{remram}, Proposition \ref{nScholzSol} we have the next result.

\begin{prop}\label{ScholzSol}
Suppose $\zeta_l \notin K$ and $K$ has no ideal classes of order $l^2$. Given a central embedding problem $(G_K, \rho, \alpha)$ with $G$ an $l$-group, cyclic $C$ and $\rho$ $l^N$-Scholz. If the extension is split or of Frattini type, then there exists an $l^N$-Scholz solution $\varphi$ and a prime $\q$ of $K$ such that 
\[Ram(\varphi) \cup T  = Ram(\rho) \cup T \cup \{\q\}. \]
\end{prop}

Define the lower central series $\{G_i\}$ of $G$ by $G_1 = G$, $G_{i+1}:=$ the commutator subgroup $[G_i,G]$,  $i\geq 1$. If $G$ is nilpotent, the smallest positive integer $c$ such that $G_{c+1} = \{1\}$ is called the nilpotency class of $G$. Our main result below generalizes Proposition 2.5 of \cite{Plans} who considers only the case $K=\Q$ and improves the result of Theorem 7.4 of \cite{GJ} when the kernel $C$ of the embedding problem is not of prime order. 

\begin{thm} \label{Ramlbound}
Given a number field $K$, a prime $l$, and an $l$-group $G$ of nilpotency class $c$. If $G$ is nonabelian, suppose $\zeta_l \notin K$ and $K$ has no ideal classes of order $l^2$. Then 
\[ minram_K(G) \leq d(G) + |T| + \sum_{i=2}^{c-1}d(G_i/G_{i+1}). \]
\end{thm}

\begin{rem}

1. This bound may be achieved by a tamely ramified extension $L/K$ with $G \cong Gal(L/K)$. 
2. If $G$ is of nilpotency class 2, 

\[minram_K(G) \leq d(G) + |T|. \]
3. If we allow $K$ to have ideal classes of order $l^2$, then the bound has the form of \cite{GJ}
\[minram_K(G) \leq g + |T|, \hs |G|=l^g.\]
\end{rem}

\proof

As in Proposition 2.5 of \cite{Plans} we use induction on $i$ for a central embedding problem 
\[\es {G_i/G_{i+1}} {G/G_{i+1}} {G/G_{i}} . \]

For $i=1$, by Proposition \ref{SplitCase} the embedding problem has an $l^N$-Scholz solution  with at most $d(G^{ab})=d(G)$ ramified primes. For $i \geq 1
$, each extension is of Frattini type, and we may break the $i$-th problem up into $d(G_i/G_{i+1})$ cyclic Frattini problems. As 
shown in Proposition \ref{ScholzSol}, each such problem may be solved at the cost of one more ramified prime. And since we can make the solution $l^N$-Scholz at each stage, it is guaranteed that we may solve the next embedding problem. \done

\section{Ramification bound on nilpotent groups \label{NilpotentSection}} 

We use the notation that $a$ is the product of the primes dividing the order of $G$ and integer $N$ satisfies $a^N$ is a multiple of the exponent of $G$. The purpose of this section is to extend Theorem \ref{Ramlbound} to groups  $G = \prod_l G_l$ that are the direct product of their Sylow $l$-subgroups $G_l$, that is {\emph {nilpotent groups}}. Assume $\zeta_l \notin K$ for all $l$ dividing $|G|$. We will obtain $G$ by a sequence of
 central embedding extensions with cyclic kernel; each of these extensions is a "product" of central extensions of $l$-groups as in sections $6$ and $7$. 
The nilpotent case was initially handled in the first author's thesis \cite{MarkinThesis}. In this section we obtain an improved bound on $minram_K(G)$ for fields $K$ which do not contain ideal classes of order $l^2$, where $l \mid |G|$.

The first step is to define a set $T$ (as small as possible) of primes of $K$ that contains an $l^N$-exceptional  set $T_l$ of primes for each $l$ dividing $|G|$. 

Let  $$\Omega_l = K(\sqrt[l^N]{E}, \sqrt[l]{V(l)}) $$ as in \ref{omega} and let $\hat \Omega = \prod_{l \mid a}\Omega_l$. Since
$Gal(\Omegae(\mu_{a^N})/K(\mu_{a^N}))$ is an $l$-group, we have  

\begin{equation}
Gal(\hat \Omega /K(\mu_{a^N})) \cong \prod_{l | a} Gal(\Omega_l (\mu_{a^N})/K(\mu_{a^N})). \label{OmegaIso}\end{equation}

Using the isomorphism of \ref{OmegaIso} we define 

$$\sigma_i = \prod_{l \mid a}\sigma_i(l), \hs  1 \leq i \leq s \hs \rm{and} $$
$$\tau_j = \prod_{l \mid a}\tau_j(l), \hs  1 \leq j \leq r \hs $$
elements of $Gal(\hat \Omega/ K(\mu_{a^N}))$. 
Here $r = \max_{l \mid a}r_l$ and we set $\tau_j(l)=1$ if $r_l < j \leq r$. 
By Chebotarev's theorem, in $K$ there is a set of $s+r$ prime ideals $T = \{ \p_i, \q_j: 1\leq i \leq s, 1\leq j\leq r \}$ disjoint from any given finite set such that 

$$\Frob {\hat \Omega / K}{\p_i} = C({Gal(\hat \Omega/K)},\sigma_i), \hs 1\leq i \leq s \hs {\rm{and}} $$ 
$$ \Frob {\hat \Omega / K}{\q_j} = C({ Gal(\hat \Omega/K)},\tau_j), \hs 1 \leq j \leq r.$$

Here $C({ Gal(\hat \Omega/K)},\gamma)$ denotes the conjugacy class  of $\gamma$ in $ Gal(\hat \Omega/K)$. 
By the properties of the Frobenius, the restriction to $\Omega_l$ of $\sigma_i$ (resp. $\tau_j$) is $\sigma_i(l)$ (resp. $\tau_j(l)$) for each $l$ dividing $a$.

\begin{lem}\label{nonsplitsol}
We continue the notation of Corollary \ref{liftinglemma} and  Lemma \ref{FieldIsoLemma}. For each $l | a$, let $L_l$ be an $l^N$-Scholz
$l$-extension of $K$ fixed by the kernel of homomorphism
$\rho_l: G_K \rightarrow \bar G_l$ and let $(G_K,\rho_l,
\alpha_l)$ be a Frattini central embedding problem as in
(\ref{cep}). Assume for all $l \mid a$, that $\zeta_{l}$ is not in $ K $ and
the exponent of $G_l$ divides $l^N$. When $|ker(\alpha_l)|>l$, assume additionally that no ideal class of $K$ has order $l^2$. Then for each $l | a$ there exists
a solution \[ \phi_l: G_K \rightarrow G_l \] for which $Ram(\phi_l)
\subseteq Ram(\rho_l)\cup T$.
\end{lem}

\proof The existence of any solution is Lemma \ref{ExistLemma}. Our set of primes $T$ contains $l^N$-exceptional subsets $T_l$, hence we may apply Theorem \ref{remram} to get a solution $\phi_l$ such that  $Ram(\phi_l) \subseteq
Ram(\rho_l)\cup T$ for all primes $l \mid a$.  \done

In the next lemma we apply Corollary \ref{corlift} to find a
single prime $\q$ that we use to lift local characters indexed by
$l | a$.

\begin{lem}

\label{mainlemma} Let $S$ be a finite set of primes of $K$ that
contains $S_0$. For each prime $l | a$, we are given integers $e_l$, $N \geq e_l \geq 1$, Galois $l$-extension $L_l/K$, character $\chi_{v,l} : K_v^\times \rightarrow \mu_{l^{e_l}}$ for all
$v \in S$. 
Assume that $K$ does not contain $\zeta_l$ for each $l | a$. There exists a prime ideal $\q$ of $K$ outside $S$ and idele class characters

$\chi_l: J_K/{K^\times} \rightarrow \mu_{l^{e_l}}$ such that conditions (1)-(4) hold for all $l \mid a$:

\begin{itemize}
\item $\q$ splits completely in $L_l(\mu_{l^{N}})/K$
\item
$\chi_l|_{\Kvs}=\chi_{v,l}$ for all $v \in S$. 
\item

$\chi_l(U_\q) = \mu_{l^{e_l}}$.

\item $\chi_l(U_v) = 1 $ for all $ v \notin S \cup \{\q\}$. 
\end{itemize}
\end{lem}

\proof
 Let $\N_l$ denote the field $ L_l(\rKS)$, $\N = \prod_{l | a}\N_l$ , $\Gamma_l = Gal( \N_l/K)$ and $\Gamma = Gal( \N/K)$.

In Lemma \ref{mScholzSoln}, for all $l \mid a$ we have defined a special prime $\q_l$ (not to be confused with  $\q_i$'s defined in Definition \ref{ExcPrimes}). Define $\sigma_l \in \Gamma_l $ by $\Frob{\N_l/K}{\q_l} = C(\Gamma_l, \sigma_l)$.
Next we show that a single prime $\q$ can be chosen.  By Lemma \ref{corlift} there exists an element ${\sigma} \in \Gamma$
whose restriction to ${\N_l}$ equals $\sigma_l$ for all $l | a$.
By Chebotarev's  theorem, there exists a prime $\q$ of $K$ 
outside $S$ such that $\Frob{\N/K}{\q} = C(\Gamma, \sigma)$. By restriction $\Frob{\N_l/K}{\q} = C(\Gamma_l, \sigma_l)$ for all $l \mid a$ and conditions (1)-(4) of \ref{mainlemma} are satisfied.
\done
        
\begin{rem}
The method by which we replaced $\{\q_l: l \mid a\}$ by $\q$ is similar to that where we replaced $\{ T_l: l\mid a\}$ by $T$.  
\end{rem}

\begin{thm}
\label{NilThm}

Given a number field $K$ and a finite nilpotent group $G$ of class c. If $G$ is nonabelian, suppose 
$\gcd(|G|,|\mu_K|)=1$ and assume for all primes dividing $|G|$ that the ideal class group of $K$ has no elements of order $l^2$. Then 

\[minram_K(G) \leq d(G) + (r+s) + \sum_{i=2}^{c-1} d(G_i/G_{i+1}).
\]

Here $s = \Z$-rank of units of $K$ and $r=\max_{l\mid |G|} \{dim Cl(K)_l\}$. 
\end{thm}

\proof

By Corollary \ref{AbG} it remains to prove the result for nonabelian groups $G$. Since $G$ is nilpotent, for each $l$ dividing $|G|$, we may apply Propositions \ref{nScholzSol}, \ref{ScholzSol}, Theorem \ref{Ramlbound} inductively. By Lemma \ref{mainlemma} there exists a single prime $\q$ to which Proposition \ref{nScholzSol} may be applied, and the conclusion follows. \done


\section{Schur Extensions}
In this section we use Fr{\"o}hlich's result on realizing the Schur multiplicator without additional ramification to verify Boston's conjecture \cite{BM} for a certain class of $l$-groups given by central extensions
$$ \es {\M(\Gamma)} G {\Gamma}.$$ 
In addition Theorem \ref{SGl3} confirms Boston's conjecture for a particular $G$ of exponent $l$, and includes the determination of the central class field of any finite abelian extension $L/\Q$ of exponent $l$ that is its own genus field. 

The group $\M(\Gamma)$ is the Schur multiplicator of a profinite group $\Gamma$ as defined in \cite{Fro}. 

\begin{defn}
Suppose $M \supseteq L \supseteq K$ are number fields with $M/K$ and $L/K$ Galois extensions. Let $M'$ be the maximal central extension of $L/K$ in $M$ and let $E$ be the maximal abelian extension of $K$ in $M$. Fr{\"o}hlich defines a certain surjective homomorphism

\begin{equation}
\M(Gal(L/K)) \rightarrow Gal(M'/EL).  \label{NrM}\end{equation}

If it is an isomorphism, one says that $M$ \emph{realizes the multiplicator $\M(Gal(L/K))$}. 
\end{defn}

\begin{rem}

If central extensions $M_1$ and $M_2$ for $L/K$ both realize the multiplicator of $Gal(L/K)$, the Galois groups $Gal(M_i/K)$, $i=1,2$ need not be isomorphic. 

\end{rem}

Fr{\"o}hlich in Proposition 3.2 of \cite{Fro} proves that if $L/K$ is an extension of finite degree, there is a finite degree central extension $M$ of $L/K$ that realizes $\M(Gal(L/K))$. 

For a prime $l$ and a finite set of primes $S$ of $K$, $K(l,S)$ denotes the maximal $l$-extension field of $K$ with ramification restricted to $S$, and $K(l,S)^{ab}$ is the maximal abelian subextension of $K(l,S)$. If $S$ contains no divisors of $l$, then the degree $[K(l,S)^{ab}:K]$ is finite. 
From now on suppose $L/K$ is a finite degree $l$-extension, so $\M(Gal(L/K))$ is a finite abelian $l$-group. Let $S$ be the set of primes of $K$ ramified in $L$. For $K=\Q$ (resp. $K$  imaginary quadratic  with $\zeta_l \notin K$), Fr{\"o}hlich in Corollary 2 of Theorem 3.13 in \cite{Fro}, (resp. Watt in  Theorem 3.1 of \cite{Watt}) proved in addition there exists such an extension $M$  that is ramified at worst at primes above $S$ (the key result on which Theorem \ref{FrW} is based). Since $M$ is central for $L/K$, we have $M=M'$. Furthermore if $L \supseteq K(l,S)^{ab}$, then $L \supseteq E$, so $EL=L$ and \ref{NrM} asserts 

\[\M(Gal(L/K)) \cong Gal(M/L). \]

\begin{rem}
Fr{\"o}hlich does not require $L/K$ to be an $l$-extension. 
\end{rem}

Thus from the results of Fr{\"o}hlich and Watt we have the following theorem. 

\begin{thm}
Let $K$ be $\Q$ or imaginary quadratic  with  $\zeta_l \notin K$ and let $L/K$ be a finite Galois $l$-extension tamely ramified only at $S$; suppose $L\supseteq K(l,S)^{ab}$. Then there exists a central extension $M$ of $L/K$  with $Ram(M/K) \subseteq S$ such that $\M(Gal(L/K)) \cong Gal(M/L)$. 
\label{FrW}\end{thm}
\done

\begin{rem}
We may apply the theorem repeatedly by replacing the extension $L/K$ by $M/K$. 
\end{rem}

\begin{rem}
Since the number of generators of $Gal(K(l,S)/K)$ equals the number of generators of $Gal(K(l,S)^{ab}/K)$, we have that $Gal(L/K)$ and $Gal(M/K)$ have the same number of generators. 
\end{rem}

\begin{thm}
Let $l$ be an odd prime and let $G$ be a group generated by $x_1, \ldots , x_n$ subject to relations $x_i^l=1$, $[x_i,x_j] \in Z(G)$ = center of $G$ for all $i,j$. There exists a tamely ramified extension $M/\Q$ such that $Gal(M/\Q) \cong G$ and $|Ram(M/\Q)|=n$. Moreover $M$ is the central class field of its maximal abelian subfield $M^{ab}$.  

\label{SGl3}
\end{thm}

\proof

Note that $\M(G/Z(G)) \cong Z(G)$.

By Proposition 2.5 of \cite{Plans} there exists an extension $M/\Q$ such that $Gal(M/\Q) \cong G$ which is (tamely) ramified at $n$ primes. Moreover we have $Ram(M/\Q)=Ram(M^{ab}/\Q)$. 

Next we show that $M/M^{ab}$ is actually unramified.
Suppose a prime ramifies in $M/M^{ab}$. As it also ramifies in $M^{ab}/\Q$, it follows that its ramification index is divisible by $l^2$. Note that since $M/\Q$ is tamely ramified, all ramification groups are cyclic. But that gives a contradiction to the exponent $l$ of $G$, thus $M$ is unramified over $M^{ab}$.

Clearly $Gal(M^{ab}/\Q)$ is the direct product of its inertia subgroups; by a standard result $M^{ab}$ is its own genus field. As $Gal(M/{M^{ab}})$ lies in the center of $Gal(M/\Q)$, $M$ is contained in the central class field $Z$ of $M^{ab}$. From e.g. \cite{Horie}, the degree of $Z$ over $M^{ab}$ divides $l^{n(n-1)/2} = 
[M:M^{ab}]$. Hence $M=Z$. 

\done

\begin{rem} In the theorem above $M$ realizes the multiplicator $\M(Gal(M^{ab}/\Q))$. 

\end{rem}
\begin{rem}
Fr{\"o}hlich \cite{Fro} Theorem 5.1 has determined the Galois group of the central class field of any abelian $l$-extension of $\Q$ by a different method. 
\end{rem}
The second author thanks Marcin Mazur for a useful discussion regarding Lemma \ref{KL}. 

\bibliographystyle{plain}
\bibliography{bib}

\end{document}